\documentclass[12pt,centertags]{amsart}
\usepackage{amsmath,amstext,amsthm,a4,amssymb,amscd}
\usepackage[mathscr]{eucal}
\usepackage{mathrsfs}
\usepackage{epsf}
\numberwithin{equation}{section}

\newtheorem{thm}{Theorem}[section]

\theoremstyle{definition}

\theoremstyle{definition}

\theoremstyle{definition}
\newtheorem{defn}[thm]{Definition}
\newcommand{\be}{\begin{eqnarray}}
\newcommand{\ee}{\end{eqnarray}}

\newcommand{\comment}[1]{}

\begin{document}

\title{Positive scalar curvature and connected sums}
\author{Guangxiang Su \ and \ Weiping Zhang}

\address{Chern Institute of Mathematics \& LPMC, Nankai
University, Tianjin 300071, P.R. China}
\email{guangxiangsu@nankai.edu.cn,\ weiping@nankai.edu.cn}

\begin{abstract}  Let $N$ be a closed enlargeable manifold in the sense of    Gromov-Lawson and $M$ a   closed spin manifold of equal dimension,  a famous theorem of Gromov-Lawson states that  the connected sum $M\# N$   admits no   metric of positive scalar curvature. 
We present a potential generalization of this result to the case where $M$ is nonspin. We use index theory for Dirac operators to prove our result. 
￼
\end{abstract}

\maketitle

\setcounter{section}{-1}

\section{Introduction} \label{s0}

It has been an important subject in differential  geometry  to study when a smooth manifold carries a Riemannian metric of positive scalar curvature (cf. \cite[Chap. IV]{LaMi89}).  A famous theorem  of Gromov and Lawson \cite{GL80}, \cite{GL83} states that an enlargeable   manifold (in the sense of \cite{GL83})  does not carry a metric of positive scalar curvature.

\begin{defn}\label{t0.1} (Gromov-Lawson  \cite[Definition 5.5]{GL83})
One calls a closed manifold $N$ (carrying a metric $g^{TN}$)  an enlargeable manifold if for any $\epsilon>0$, there is a   covering manifold $\widehat N_\epsilon\rightarrow N$, with $\widehat N_\epsilon$ being spin,  and a smooth map $f:\widehat N_\epsilon\rightarrow S^{\dim N}(1)$ (the standard unit sphere), which is constant near infinity and has non-zero degree, such that for any $X\in \Gamma(T\widehat N_\epsilon)$, $|f_*(X)|\leq \epsilon |X|$.
\end{defn}

It is clear that the enlargeability does not depend on the metric $g^{TN}$.

We assume from now on that $N$ is a  closed enlargeable manifold. Let $M$ be a closed    manifold such that there is a closed codimension two submanifold $W\subset M$ such that $M\setminus W$ is spin. Without loss of generality, we assume that $\dim M=\dim N=n$ is even.  Let   $h^{TN}$  be a metric on   $TN$.

We fix a point $p\in N$. For any $r\geq 0$, let $B^N_p(r)=\{y\in N\,:\, d(p,y)\leq r\}$. Let $a_0>0$ be a fixed sufficiently small number. Then the connected sum $M \# N$ can be constructed so that  the hypersurface $\partial B^N_p(a_0)$, which is the boundary of $B^N_p(a_0)$,  cuts $M \# N$ into two parts: the part $N\setminus B^N_p(a_0)$ and the rest part coming from $M  $ (by attaching the boundary of a ball in $M\setminus W$ to $\partial B^N_p(a_0)$). 

Let $\varphi: M \# N\rightarrow [0,1]$ be an arbitrary smooth  function  such that $\varphi\equiv 1$ on $N\setminus  B^N_p(a_0)$ and ${\rm Supp}(\varphi)\subseteq M\# N\setminus W$. 
The main result of this paper can be stated as follows.

\begin{thm}\label{t0.2}  There is no metric $g^{T(M \# N)}$ on $T(M \# N)$ such that the associated scalar curvature $k^{T(M \# N)}$ verifies the following inequality on ${\rm Supp}(\varphi)$,
\begin{align}\label{0.1}
c-6|d\varphi|^2_{ g^{T(M \# N)}}\geq {\rm max}\left\{0, \frac{3c}{2}-\frac{k^{T(M \# N)}}{4}\right\}
\end{align} 
for some constant $c>0$. 
\end{thm}

When   $M$ is spin, one can take $W=\emptyset$, $\varphi\equiv 1$ on $M\# N$ and   $c>0$ small enough to  recover the  theorem of Gromov-Lawson \cite{GL80}, \cite{GL83}  mentioned at the begining.

 Our proof of Theorem \ref{t0.2}   is index theoretic and   is inspired by    \cite{Z17a}, where a new proof of the above mentioned Gromov-Lawson theorem is given without using index theorems on noncompact manifolds. 
  The details  will be carried out in Section \ref{s2}.

\section{Proof of Theorem \ref{t0.2}}\label{s2}

Assume   there is a metric  $g^{T(M \# N)}$  on $T(M \# N)$ such that  (\ref{0.1}) holds for $c=\alpha^2>0$.


 For any $\epsilon>0$, let $\pi:\widehat N_\epsilon\rightarrow N$ be a  covering manifold verifying Definition \ref{t0.1}, carrying lifted geometric data from that of $N$. Let $a_0>0$ be small enough so that for any $p',\,q'\in\pi^{-1}(p)$ with $p'\neq q'$, $\overline{B^{\widehat N_\epsilon}_{p'}(4a_0)}\cap \overline{B^{\widehat N_\epsilon}_{q'}(4a_0)}=\emptyset$.\footnote{Here and in what follows, the involved balls are determined by $h^{TN}$.}   It is clear that one can choose $a_0>0$ not depending on $\epsilon$. 

Let $h: N\rightarrow N$ be a smooth map such that $h={\rm Id}$ on $N\setminus B^N_p(3a_0)$, while $h( B_p^N(2a_0))=\{p\}$.  
It lifts to a map $\widehat h:\widehat N_\epsilon\rightarrow \widehat N_\epsilon$ verifying  that $\widehat h={\rm Id}$ on $\widehat N_\epsilon\setminus \bigcup_{p'\in\pi^{-1}(p)}B_{p'}^{\widehat N_\epsilon}(3a_0)$, while for any $p'\in\pi^{-1}(p)$, $\widehat h(B_{p'}^{\widehat N_\epsilon}(2a_0))=\{p'\}$.

Let $f:  \widehat N_\epsilon\rightarrow S^{n}(1)$ be as in Definition \ref{t0.1}. Set $\widehat f=f\circ \widehat h:\widehat{  N}_\epsilon\rightarrow S^n(1)$. Then ${\rm deg}(\widehat f)={\rm deg}(f)$ and  there is a   constant $c'>0$ such that for any $X\in\Gamma(T\widehat {  N}_\epsilon)$, one has 
\begin{align}\label{1.1}
\left|\widehat f_*(X)\right|\leq c'\,\epsilon\, |X|.
\end{align} 


To simplify the presentation, we assume first that each $\widehat N_\epsilon$ is compact, i.e., $N$ is a {\it compactly} enlargeable manifold. 

  Since $M\setminus W$ is spin,   one can construct  a compact spin manifold with boundary $M_W\subset M\setminus W$ such that $\partial M_W \subseteq M\setminus {\rm Supp}(\varphi)$. Let $M'_W$ be another copy of $M_W$. Then one gets a closed spin manifold by gluing $M_W$ and $M'_W$ along the boundary. We denote  the resulting double by  $\widetilde M_W$. 
Then one can extend the connected sum $M_W\# N$ to $\widetilde M_W\# N$ obviously.  It lifts natually to $\widehat N_\epsilon$ where near each $p'\in \pi^{-1}(p)$, we do the  lifted connected sum. We denote the resulting manifold by $\widehat M_W\# \widehat N_{\epsilon}$. 
Clearly, the metric $g^{T(M \# N)}$ induces a metric $g^{T(\widetilde M_W \#   N) }$ such that $g^{T(\widetilde M_W \#   N )}|_{M_W \#   N }= g^{T(M \# N)}|_{M_W \#   N } $. Let $k^{T(\widetilde M_W \#   N )}$ be the associated scalar curvature. They determine the corresponding metric $g^{T(\widehat M_W \#\widehat   N_\epsilon) }$ and scalar curvature $k^{T(\widehat M_W \#\widehat   N_\epsilon) }$ on $\widehat M_W \#\widehat   N_\epsilon $. 

The cut-off function    $\varphi$   extends   to $\widetilde  M_W\#  N $ by setting    $\varphi(M'_W)=0$. 
It lifts  to $\widehat M_W\# \widehat N_{\epsilon}$ obviously and we still denote the lifting by $\varphi$.

We extend $\widehat f:\widehat  N_{\epsilon}\rightarrow S^{n}(1)$ to $\widehat f:\widehat M_W\# \widehat N_{\epsilon}\rightarrow S^{n}(1)$ by setting $\widehat f(\widetilde  M_W\#  B_{p'}(4a_0))=f(p')$ for any $p'\in\pi^{-1}(p)$. 


Following \cite{GL80}, \cite{GL83},  
let $(E_0,g^{E_0} )$ be a Hermitian vector bundle on $S^{n}(1)$   carrying a Hermitian connection $\nabla^{E_0}$ such that
\begin{align}\label{1.6}
\left\langle {\rm ch}\left(E_0\right),\left[S^n(1)\right]\right\rangle\neq 0.
\end{align}
  Let $(E_1={\bf C}^k|_{S^{n}(1)},g^{E_1},\nabla^{E_1})$, with $k={\rm rk}(E_0)$, be the canonical Hermitian trivial vector bundle on $S^{n}(1)$.  

For any $p'\in\pi^{-1}(p)$, let $v_{f(p')} :\Gamma(E_0|_{f(p')})\rightarrow \Gamma (E_1|_{f(p')})$ be    an isometry.  Let $v_{f(p')}^*: \Gamma(E_1|_{f(p')})\rightarrow \Gamma(E_0|_{f(p')})$ be the adjoint of $v_{f(p')}$ with respect to   $g^{E_0}|_{f(p')}$ and $g^{E_1}|_{f(p')}$. Set 
\begin{align}\label{1.11}
V_{f(p')}=v_{f(p')}+v^*_{f(p')}. 
\end{align}

Let $(\xi,g^\xi,\nabla^\xi)=(\xi_0\oplus\xi_1,g^{\xi_0}\oplus g^{\xi_1},\nabla^{\xi_0}\oplus\nabla^{\xi_1})=(\widehat f ^*E_0\oplus \widehat f^*E_1, \widehat f^*g^{E_0}\oplus \widehat f^*g^{E_1},\widehat f^*\nabla^{E_0}\oplus \widehat f^*\nabla^{E_1})$ be the  ${\bf Z}_2$-graded Hermitian vector bundle with Hermitian connection over $\widehat M_W\# \widehat  N_{\epsilon}$. Let $R^\xi=(\nabla^\xi)^2$ be the curvature of $\nabla^\xi$.

\comment{
Set $V_{\widehat f}=\widehat f^*V$. Then  
\begin{align}\label{1.10}
 V_{\widehat f} ^2=1  \ \ \ {\rm and}\ \ \ \left[\nabla^\xi,V_{\widehat f}\right]=0
\end{align}
on $\widehat  N_{\epsilon,H}'\cup (H\times[0,\alpha_0))\cup \bigcup_{p'\in\pi^{-1}(p)\cap K}B_{p'}^{\widehat  N_\epsilon}(2a_0)$. 

}

Let $D^\xi:\Gamma(S(T(\widehat M_W\# \widehat  N_{\epsilon}))\widehat\otimes\xi)\rightarrow \Gamma(S(T(\widehat M_W\# \widehat  N_{\epsilon}))\widehat\otimes \xi)$ be the canonical  (twisted) Dirac operator (cf. \cite{LaMi89}) associated to $(T(\widehat M_W\# \widehat  N_{\epsilon}),g^{T(\widehat M_W\# \widehat  N_{\epsilon})})$ and $(\xi,g^\xi,\nabla^\xi)$. Let $D^\xi_\pm: \Gamma((S(T(\widehat M_W\# \widehat  N_{\epsilon}))\widehat\otimes\xi)_\pm)\rightarrow \Gamma((S(T(\widehat M_W\# \widehat  N_{\epsilon}))\widehat\otimes \xi)_\mp)$ be the obvious restrictions, where $(S(T(\widehat M_W\# \widehat  N_{\epsilon}))\widehat\otimes\xi)_+=S_+(T(\widehat M_W\# \widehat  N_{\epsilon}))\otimes\xi_0\oplus S_-(T(\widehat M_W\# \widehat  N_{\epsilon}))\otimes \xi_1$, while $(S(T(\widehat M_W\# \widehat  N_{\epsilon}))\widehat\otimes\xi)_-=S_-(T(\widehat M_W\# \widehat  N_{\epsilon}))\otimes\xi_0\oplus S_+(T(\widehat M_W\# \widehat  N_{\epsilon}))\otimes \xi_1$. By the Atiyah-Singer index theorem \cite{ASI} (cf. \cite{LaMi89})  and \cite{GL83}, one has 
\begin{multline}\label{1.12}
{\rm ind}\left(D^\xi_+\right)=\left\langle\widehat A\left(T\left(\widehat M_W\# \widehat  N_{\epsilon}\right)\right)\left({\rm ch}\left(\xi_0\right)-{\rm ch}\left(\xi_1\right)\right),\left[\widehat M_W\# \widehat  N_{\epsilon}\right]\right\rangle 
\\
= ({\rm deg}(f))\left\langle {\rm ch}\left(E_0\right),\left[S^{n}(1)\right]\right\rangle. 
\end{multline}
 
 Following \cite[p. 115]{BL91}, let $\varphi_1,\,\varphi_2:\widehat M_W\# \widehat  N_{\epsilon}\rightarrow [0,1]$ be defined by 
\begin{align}\label{1.14}
\varphi_1=\frac{\varphi} {\left(\varphi^2+(1-\varphi)^2\right)^{\frac{1}{2}}},\ \ \ \varphi_2=\frac{1-\varphi} {\left(\varphi^2+(1-\varphi)^2\right)^{\frac{1}{2}}}.
\end{align}
Then $\varphi_1^2+\varphi_2^2=1$.

Recall that  $\alpha^2=c>0$. Let $D_\alpha^\xi: \Gamma(S(T(\widehat M_W\# \widehat  N_{\epsilon}))\widehat\otimes\xi)\rightarrow \Gamma(S(T(\widehat M_W\# \widehat  N_{\epsilon}))\widehat\otimes \xi)$ be the deformed operator defined by
\begin{align}\label{1.13}
D^\xi_\alpha=D^\xi + \alpha\sum_{p'\in \pi^{-1}(p)} \varphi_2 {V_{ p' }} ,
\end{align}
where  $V_{p'} = \widehat f^{-1}(V_{f(p')})$ lives on the lift of $M_W\# B_{p}(4a_0)$ near $p'$. 

 From (\ref{1.13}), it is easy to verify that 
\begin{align}\label{1.1x}\left(D_\alpha^\xi \right)^2=\left(D^\xi\right)^2 +\alpha\sum_{p'\in \pi^{-1}(p)} c\left(d\varphi_2\right)V_{p'}+\alpha^2\varphi ^2_2 .
\end{align}
Thus for any 
 $s\in \Gamma(S(T(\widehat M_W\#\widehat  N_{\epsilon}))\widehat\otimes\xi)$, we have
\begin{multline}\label{1.7}
\left\| D_\alpha^\xi s\right\|^2=\|D^\xi s\|^2+\alpha\sum_{p'\in \pi^{-1}(p)}\left\langle c\left(d\varphi_2\right)V_{p'} s,s\right\rangle+\alpha^2  \left\|\varphi_2 s\right\|^2\\
=\left\|\varphi_1 D^\xi s\right\|^2+\left\|\varphi_2 D^\xi s\right\|^2+\alpha\sum_{p'\in \pi^{-1}(p)}\left\langle c\left(d\varphi_2\right)V_{p'} s,s\right\rangle+  \alpha^2\left \|\varphi_2 s\right\|^2.
\end{multline}

By direct computation we have for $i=1,\, 2$ that 
\begin{align}\label{1.8}
\left\|\varphi_i D^\xi s\right\|^2= \left\|D^\xi\left(\varphi_i s\right)\right\|^2 -\left\|\left|d\varphi_i\right|s\right\|^2-\left\langle D^\xi s, {c\left(d\varphi_i^2\right)\over 2} s\right\rangle-\left\langle {c\left(d\varphi_i^2\right)\over 2} s,D^\xi s\right\rangle.
\end{align}

By (\ref{1.7}) and (\ref{1.8}), we have
\begin{align}\label{1.2}
\left\| D^\xi_\alpha s\right\|^2=\sum_{i=1}^2\left(\|D^\xi\left(\varphi_i s\right)\|^2  -\left\|\left|d\varphi_i\right|s\right\|^2\right) 
+\alpha \sum_{p'\in \pi^{-1}(p)} \left\langle c\left(d\varphi_2\right)V_{p'} s,s\right\rangle+\alpha^2 \left\|\varphi_2s\right\|^2.
\end{align}

By (\ref{1.14}), we have
\begin{align}\label{1.1a}
d\varphi_1= {\varphi_2 d\varphi \over {\varphi^2 +(1-\varphi)^2}},\ \ \ \ d\varphi_2= -{\varphi_1 d\varphi \over {\varphi^2 +(1-\varphi)^2}}.
\end{align}

Let $e_1,\,\cdots,\,e_n$ be an orthonormal basis of $(T(\widehat M_W\#\widehat N_\epsilon),g^{T(\widehat M_W\#\widehat N_\epsilon)})$. 
By (\ref{1.1}), (\ref{1.2}), (\ref{1.1a}),  the Lichnerowicz formula \cite{L63} (cf. \cite{LaMi89}) and proceed as in \cite{GL83}, one has
\begin{multline}\label{1.1b}
\left\| D^\xi_\alpha s\right\|^2=-\left\langle \Delta\left(\varphi_1 s\right),\varphi_1 s\right\rangle+\left\langle {k^{T (\widehat M_W \#\widehat   N_\epsilon ) }\over 4}\varphi_1 s,\varphi_1 s\right\rangle+\left\|D^\xi\left(\varphi_2 s\right)\right\|^2
\\
-\left\|{|d\varphi|\over{\varphi^2+(1-\varphi)^2}}s\right\|^2-\alpha\sum_{p'\in \pi^{-1}(p)}
\left\langle{\varphi_1c(d\varphi)V_{p'} \over{\varphi^2+(1-\varphi)^2}}s,s\right\rangle+\alpha^2 \left \|\varphi_2s\right\|^2
\\
+{1\over 2}\sum_{i,\,j=1}^n \left\langle c\left(e_i\right)c\left(e_j\right)R^\xi\left(e_i,e_j\right)s,s\right\rangle
\\
\geq \left\langle\left( {k^{T(\widehat M_W \#\widehat   N_\epsilon) }\over 4}\varphi_1^2+\alpha^2\varphi^2_2  -{|d\varphi|^2\over{(\varphi^2+(1-\varphi)^2)^2}}-\alpha\sum_{p'\in \pi^{-1}(p)} {{\varphi_1 c(d\varphi)V_{p'} }\over{\varphi^2+(1-\varphi)^2}}\right)s,s\right\rangle
\\
+{1\over 2}\sum_{i,\,j=1}^n \left\langle c\left(e_i\right)c\left(e_j\right) \widehat f^*\left( \left(\nabla^{E_0}\right)^2\left(\widehat f_*\left(e_i\right),\widehat f_*\left(e_j\right) \right)\right)s,s\right\rangle
\\
\geq\left\langle\left( {k^{T(\widehat M_W \#\widehat   N_\epsilon) }\over 4}\varphi_1^2+\alpha^2\varphi^2_2  -\left.{{\alpha^2 \varphi_1^2 }\over 2}\right|_{{\rm Supp}(d\varphi)}-{3\over 2}{{|d\varphi|^2}\over{(\varphi^2+(1-\varphi)^2)^2}}\right)s,s\right\rangle 
\\
+\left\langle O\left(\epsilon^2\right)s,s\right\rangle_{{\widehat N_\epsilon}\setminus \cup_{p'\in\pi^{-1}(p)}B_{p'}^{\widehat N_\epsilon}(a_0)}
\\
\geq\left\langle\left( {k^{T(\widehat M_W \#\widehat   N_\epsilon) }\over 4}\varphi_1^2+\alpha^2\varphi^2_2  -\left.{{\alpha^2 \varphi_1^2 }\over 2}\right|_{{\rm Supp}(d\varphi)}-{6}{{|d\varphi|^2} }\right)s,s\right\rangle
\\
+\left\langle O\left(\epsilon^2\right)s,s\right\rangle_{{\widehat N_\epsilon}\setminus \cup_{p'\in\pi^{-1}(p)}B_{p'}^{\widehat N_\epsilon}(a_0)}.
\end{multline}

 For any $x\in {\rm Supp}(d\varphi)$, if ${k^{T(\widehat M_W \#\widehat   N_\epsilon) }\over 4}- {{3\alpha^2   }\over 2}\geq 0$ at $x$, then one has  at $x$ that, in view of (\ref{0.1}),
\begin{multline}\label{1.22}
 {k^{T(\widehat M_W \#\widehat   N_\epsilon) }\over 4}\varphi_1^2+\alpha^2\varphi^2_2  -\left.{{\alpha^2 \varphi_1^2 }\over 2}\right|_{{\rm Supp}(d\varphi)}-{6}{{|d\varphi|^2} }
\\
= \left({k^{T(\widehat M_W \#\widehat   N_\epsilon) }\over 4} -  {{3\alpha^2   }\over 2}\right) \varphi_1^2+\alpha^2- {6}{{|d\varphi|^2} }\geq 0,
\end{multline}
while if ${k^{T(\widehat M_W \#\widehat   N_\epsilon) }\over 4}- {{3\alpha^2   }\over 2}\leq 0$ at $x$, then by (\ref{0.1}), one has at $x$ that
\begin{multline}\label{1.23}
 {k^{T(\widehat M_W \#\widehat   N_\epsilon) }\over 4}\varphi_1^2+\alpha^2\varphi^2_2  -\left.{{\alpha^2 \varphi_1^2 }\over 2}\right|_{{\rm Supp}(d\varphi)}-{6}{{|d\varphi|^2} }
\\
= \left(    {{3\alpha^2   }\over 2}-{k^{T(\widehat M_W \#\widehat   N_\epsilon) }\over 4}\right) \varphi_2^2+{k^{T(\widetilde M_W \#   N )}\over 4}-{\alpha^2\over 2}- {6}{{|d\varphi|^2} }\geq 0. 
\end{multline}

On the other hand, (\ref{0.1}) implies that
\begin{align}\label{1.24}
 {k^{T(\widehat M_W \#\widehat   N_\epsilon) }\over 4}\geq {\alpha^2\over 2} 
\end{align}
on $\pi^{-1}(N\setminus B_{p}(a_0))={{\widehat N_\epsilon}\setminus \cup_{p'\in\pi^{-1}(p)}B_{p'}^{\widehat N_\epsilon}(a_0)}$, on which $\varphi\equiv 1$. 

From (\ref{1.1b})-(\ref{1.24}), we see that if (\ref{0.1}) holds for $c=\alpha^2$, then one has
\begin{align}\label{1.25}
 \left\| D^\xi_\alpha s\right\|^2\geq \left\langle \left({\alpha^2\over 2}+O\left(\epsilon^2\right)\right)s,s\right\rangle_{{\widehat N_\epsilon}\setminus \cup_{p'\in\pi^{-1}(p)}B_{p'}^{\widehat N_\epsilon}(a_0)},
\end{align}
which implies (when $\epsilon>0$ is small enough)
 $\ker(D_\alpha^\xi)=\{0\}$ (cf. \cite[Theorem 8.2]{BoW93}), which contradicts (\ref{1.12}) where the right hand side is nonzero.  This completes the proof of Theorem \ref{t0.2} for the case where $N$ is a compactly enlargeable manifold.

For the general case where $\widehat N_\epsilon$ is noncompact, one can combine the above arguments with the method in \cite{Z17a} to complete the proof of Theorem \ref{t0.2}. We leave the details to the interested reader.



$\ $

\noindent{\bf Acknowledgments.} This work was partially supported by NNSFC.

\def\cprime{$'$} \def\cprime{$'$}
\providecommand{\bysame}{\leavevmode\hbox to3em{\hrulefill}\thinspace}
\providecommand{\MR}{\relax\ifhmode\unskip\space\fi MR }
\providecommand{\MRhref}[2]{%
  \href{http://www.ams.org/mathscinet-getitem?mr=#1}{#2}
}
\providecommand{\href}[2]{#2}


\begin{thebibliography}{99}




\comment{

\bibitem{AH} M.~F.~Atiyah and F. Hirzebruch, {\it Spin-manifolds and
group actions}. in {\it Essays on Topology and Related Topics},
Memoires d\'edi\'es \`a Georges de Rham (ed. A. Haefliger and R.
Narasimhan), Springer-Verlag, New York-Berlin (1970), 18-28.

}


\comment{

\bibitem{APS} M.~F.~Atiyah, V.~K.~Patodi and
I.~M.~Singer, {\it Spectral asymmetry and Riemannian geometry I}.
{Proc. Camb. Philos. Soc.} {\bf 77} (1975), 43-69.

}

\bibitem{ASI} M.~F.~Atiyah  and I.~M.~Singer,
{\it The index  of elliptic operators I}. Ann. of Math. 87 (1968),
484-530.


\comment{
\bibitem{ASV} M.~F.~Atiyah  and I.~M.~Singer, {\it The index  of elliptic operators V}. Ann. of Math. 93 (1971),
139-149.

}


 \comment{

\bibitem{BeGeVe}
N.~Berline, E.~Getzler and M.~Vergne, \emph{Heat Kernels and {D}irac
 Operators}. Grundl. Math. Wiss. Band 298, Springer-Verlag, Berlin, 1992.

}

 
\comment{
  \bibitem{BC89} J.-M. Bismut and J. Cheeger, {\it $\eta$-invariants
  and their adiabatic limits.} J. Amer. Math. Soc. 2 (1989),
  33-70.

}

\bibitem{BL91}
J.-M. Bismut and G.~Lebeau, \emph{Complex immersions and {Q}uillen metrics}.
  Inst. Hautes {\'E}tudes Sci. Publ. Math. (1991), no.~74, ii+298
  pp.



\bibitem{BoW93} B.~Boo{\ss}-Bavnbek and K.~Wojciechowski,
\emph{Elliptic boundary problems for {D}irac operators},
Mathematics: Theory \& Applications, Birkh{\"a}user Boston
  Inc., Boston, MA, 1993.


\comment{

\bibitem{Bo70}
R. Bott, {\it On a topological obstruction to integrability}. in
{\it Global Analysis}. Proc. Symp. Pure Math. vol. 16, (1970),
127-131.

 \bibitem{Brav02} M.~Braverman, \emph{Index theorem for
equivariant {D}irac operators on
  non-compact manifolds}. $K$-Theory  {27} (2002), 
  61--101.

\bibitem{CCL} S. S. Chern, W. H. Chen and K. S. Lam, {\it Lectures on Differential Geometry.}  Series on University Mathematics  - Vol. 1, World Scientific, 1999. 



\bibitem{Co86}
A. Connes, \emph{Cyclic cohomology and the transverse fundamental
class of a foliation}. in {\it Geometric Methods in Operator
Algebras.} H. Araki eds., pp. 52-144, Pitman Res. Notes in Math.
Series, vol. 123, 1986.
}

\comment{

\bibitem{CoS84}
A. Connes and G. Skandalis, {\it The longitudinal index theorem for
foliations}. Publ. Res. Inst. Math. Sci. Kyoto, 20 (1984),
1139-1183.

}

\comment{

\bibitem{DZ00} X. Dai and W. Zhang, {\it Real embeddings and the Atiyah-Patodi-Singer index
theorem for Dirac operators.} Asian J. Math. 4 (2000), 775-794.

}


\comment{\bibitem{Gilkey93}
 P.~Gilkey,  \emph{On the index of geometrical operators for Riemannian
manifolds with boundary}. Adv. Math. \textbf{102} (1993), no. 2,
129--183.}


\comment{

\bibitem{Gr96} M. Gromov, {\it Positive curvature, macroscopic dimension, spectral gaps and higher
signatures}. in {\it Functional Analysis on the Eve of the 21st
Century}, Gindikin, Simon (ed.) et al.     Progress in Math. 132
(1996), 1-213, Birkh\"auser, Basel.

}

\bibitem{GL80} M. Gromov and H. B. Lawson, {\it Spin and scalar
curvature in the presence of a fundamntal group I}. Ann. of Math.
111 (1980), 209-230.

\comment{

\bibitem{GL80a} M. Gromov and H. B. Lawson, {\it The classification of simply-connected manifolds of
positive scalar curvature}. Ann. of Math. 111 (1980), 423-434.

}

\bibitem{GL83} M. Gromov and H. B. Lawson,  {\it  Positive scalar curvature and the Dirac operator on complete Riemannian manifolds}.  Publ. Math. I.H.E.S. 58 (1983), 295-408.

\comment{

\bibitem{He} S. Helgason, {\it Differential Geometry and Symmetric
Spaces}. Academic Press, 1962.




\bibitem{HR} N. Higson and J. Roe, {\it Lectures on Noncommutative
Geometry} (2000 Clay Mathematics Institute Instructional
Symposium). Slides available at
http://www.personal.psu.edu/ndh2/math/Slides.html.   

 \comment{\bibitem{HS}
M. Hilsum and G. Skandalis, {\it Morphisms $K$-orient\'es
d'espaces de feuilles et fonctorialit\'e en th\'eorie de
Kasparov.} Ann. Scient. ENS 20 (1987), 325-390.}

 \bibitem{Hi74} N. J.  Hitchin,     {\it Harmonic spinors}.
Adv. in Math. 14 (1974), 1-55.

}

\bibitem{LaMi89}
H.~B. Lawson and M.-L. Michelsohn, \emph{Spin Geometry}.   Princeton Univ. Press, Princeton, NJ, 1989.

\bibitem{L63}
A. Lichnerowicz, {\it Spineurs harmoniques.} C. R. Acad. Sci. Paris,
S\'erie A, 257 (1963), 7-9.


\comment{


\bibitem{LMZ01}   K. Liu, X. Ma and W. Zhang, {\it On elliptic genera and foliations}. Math. Res. Lett. 8 (2001), 361-376.




\bibitem{LW09}  K. Liu and Y. Wang, {\it Adiabatic limits, vanishing theorems and the noncommutative
residue}. Science in China Series A: Mathematics    52 (2009),
2699-2713.

}

 



\comment{

\bibitem{Lu71} G. Lusztig, {\it Novikov's higher signature and
families of elliptic operators}. J. Diff. Geom. 7 (1971), 229-256.

\bibitem{MZ12} X. Ma and W. Zhang, 
{\it Transversal index and $L^2$-index for manifolds with boundary}. in 
{\it Metric and Differential Geometry, The Jeff Cheeger Anniversary Volume}, pp. 299-315. Eds. X. Dai and X. Rong. Progress in Mathematics, vol. 297. Birkh\"auser Boston, Inc., Boston, MA. 2012. 

}










\comment{

\bibitem{S92} S. Stolz, {\it Simply connected manifolds of positive
scalar curvature.} Ann. of Math. 136 (1992), 511-540.
\bibitem{T}  W. P. Thurston, {\it Existence of codimension-one foliations}. Ann. of Math. 104 (1976), 249-268.


 
}




 


\comment{

\bibitem{Z04} W. Zhang, {\it Sub-signature operators, $\eta$-invariants and a Riemann-Roch theorem for flat vector bundles}. Chin. Ann. Math. 25B (2004), 7-36.

}


\bibitem{Z17a} W. Zhang,  {\it Positive scalar curvature on foliations: the enlargeability.} Preprint, arXiv:1703.04313. 

\comment{

\bibitem{MS88}  R. J. Zimmer, {\it Positive scalar curvature along the leaves}. Appendix C in 
 {\it Global Analysis
on Foliated Spaces}.  By  C. C. Moore and C. Schochet.  MSRI Publ. Vol. 9. Springer-Verlag, 1988.

}


\end{thebibliography}
\end{document}